\documentclass[12pt,a4paper]{article}
\usepackage[cp1251]{inputenc}
\usepackage[english]{babel}
\usepackage[OT1]{fontenc}
\usepackage{amsmath}
\usepackage{amsfonts}
\usepackage{amssymb}
\usepackage{xcolor}
\usepackage{url}
\usepackage{graphicx}
\usepackage{setspace} 
\newcounter{theorem}

\newcounter{theoremcounter}

\newcounter{remarkcounter}

\newtheorem{theorem}[theoremcounter]{Theorem}

\newtheorem{remark}[remarkcounter]{Remark}

\newcommand{\vp}{{\mathbf p}}
\newcommand{\vx}{{\mathbf x}}

\newcommand{\vy}{{\mathbf y}}

\newcommand{\be}{\begin{equation}}
\newcommand{\ee}{\end{equation}}
\newcommand{\ber}{\begin{eqnarray}}
\newcommand{\eer}{\end{eqnarray}}

\begin{document}

\title{On the Rate of Convergence for a Characteristic of Multidimensional Birth-Death Process}

\author{A.\,I.~Zeifman\textsuperscript{1}\and Y.\,A.~Satin\textsuperscript{2}\and K.\,M.~Kiseleva\textsuperscript{3}\and V.\,Yu.~Korolev\textsuperscript{4}}

\date{}

\maketitle

\footnotetext[1]{Department of Applied Mathematics, Vologda State
University, Vologda, Russia; IPI FRC CSC RAS, VolSC RAS;
\url{a_zeifman@mail.ru}} \footnotetext[2]{Department of
Mathematics, Vologda State University, Vologda, Russia;\\
\url{yacovi@mail.ru}} \footnotetext[3]{Department of Applied
Mathematics, Vologda State University, Vologda, Russia;
\url{ksushakiseleva@mail.ru}} \footnotetext[4]{Faculty of
Computational Mathematics and Cybernetics, Lomonosov Moscow State
University, Moscow, Russia; IPI FRC CSC RAS, Moscow, Russia;
Hangzhou Dianzi University, China; \url{vkorolev@cs.msu.ru}}

\begin{abstract}
We consider a multidimensional inhomogeneous birth-death process
(BDP) and obtain bounds on the rate of convergence for the
corresponding one-dimensional processes.
\end{abstract}

\section{Introduction}
Multidimensional birth-death processes (BDP) were objects of a
number of studies in queueing theory and other applied fields, see
\cite{akpan,ds,deng,goswami2017,goswami2019,jl,js,lee,marsan,mather,michaelides,ni,stevens}. The
problem of the product form solutions for such models was
considered, for instance, in \cite{tokb} (also, see the references
therein). If the process is inhomogeneous and transition intensities
have a more general form, then the problem of computation of any
probabilistic characteristics of the queueing model is much more
difficult.

The background of our approach is the method of investigation of
inhomogeneous BDP, see the detailed discussion and some preliminary
results in \cite{gz04,z95jap,z95spa,z06}. Estimates for the state
probabilities of one-dimensional projections of a multidimensional
BDP were studied in \cite{z15dan} and \cite{z16jms}. However, within
that methodology it was impossible to obtain estimates of the rate
of convergence, since the logarithmic norm of the operator cannot be
applied to the corresponding nonlinear systems.

\smallskip

 Here we substantially modify that approach so that it
can be used for estimation and construction of some explicit bounds
on the rate of convergence for one-dimensional projection of a
multidimensional BDP. Namely, in Section 2 we develop a simple but
efficient method for bounding the rate of convergence for an
arbitrary (may be, nonlinear, depending on the number of parameters
and so on) differential equation in the space of sequences $l_1$,
and in Section 3 we apply this method to bounding the rate of
convergence for one-dimensional projections of BDP.

\smallskip

Let ${\bf X}(t)=(X_1(t),...,X_d(t))$ be a $d$-dimensional  BDP such
that in the interval $(t,t+h)$ the following transitions are
possible with order $h$: birth of a particle of type $j$, death of a
particle of type $j$.

Let $\lambda _{j,{\bf m}}(t)$ be the corresponding birth rate (from
the state ${\bf m}=(m_1,...,m_d)=\sum_{i=1}^dm_i{\bf e}_i$ to the
state ${\bf m}+{\bf e}_j$) and $\mu _{j,{\bf m}}(t)$  be the
corresponding death intensity (from the state ${\bf
m}=(m_1,...,m_d)=\sum_{i=1}^dm_i{\bf e}_i$ to the state ${\bf
m-e}_j$). Denote $p_{{\bf m}}(t)=\Pr \left( {\bf X}(t)={\bf
m}\right)$.

\smallskip

Let now the (countable) state space of the vector process under
consideration be arranged in a special order, say $0, 1, \dots$.
Denote by $p_i(t)$ the corresponding state probabilities, and by
${\bf p}(t)$ the corresponding column vector of state probabilities.
Applying our standard approach (see details in
\cite{gz04,z95spa,z06}) we suppose in addition, that all intensities
are nonnegative functions locally integrable on $[0,\infty)$, and,
moreover, in new enumeration,
$$
\Pr\left(X( t+h) =j/X(t) =i\right) = \bigg\{q_{ij}(t)
h+\alpha_{ij}(t, h), ~~j\neq i,~~~~~
$$
\begin{equation}
1-\sum\limits_{k\neq i}q_{ik}(t) h+\alpha_{i}(t,h), ~~j=i,\bigg\}
\label{1001}
\end{equation}

\noindent where all  $\alpha_{i}(t,h)$ are $o(h)$ uniformly in
$i$, i. e. $\sup_i |\alpha_i(t,h)| = o(h)$.

We suppose that $\lambda _{j,{\bf m}}(t) \le L <\infty, ~~ \mu _{j,{\bf m}}(t)
\le M  <\infty,$ for any $j$, ${\bf m}$ and almost all $t \ge 0$.

The probabilistic dynamics of the process is represented
by the forward Kolmogorov system:
\begin{equation} \label{ur01}
\frac{d\vp}{dt}=A(t)\vp(t),
\end{equation}
\noindent where $A(t)$ is the corresponding infinitesimal
(intensity) matrix.

Throughout the paper we denote  the $l_1$-norm by $\|\cdot\|$, i. e.
$\|{\vx}\|=\sum|x_i|$, and $\|B\| = \sup_j \sum_i |b_{ij}|$ for $B =
(b_{ij})_{i,j=0}^{\infty}$.

Let $\Omega$ be the set all stochastic vectors, i. e., $l_1$-vectors
with nonnegative coordinates and unit norm. We have the inequality
$\|A(t)\|
 \le  2d\left(L+M\right)< \infty , $
\noindent for any $j$, ${\bf m}$ and almost all $t \ge 0$.  Hence,
the operator function $A(t)$ from $l_1$ into itself is bounded for
almost all $t \ge 0$ and is locally integrable on $[0;\infty)$.
Therefore we can consider (\ref{ur01}) as a differential equation in
the space $l_1$ with bounded operator.

It is well known, see \cite{dk}, that the Cauchy problem for
differential equation (\ref{ur01}) has unique solution for an
arbitrary initial condition, and  $\vp(s) \in \Omega$ implies
$\vp(t) \in \Omega$ for $t \ge s \ge 0$.

\smallskip

We recall that a Markov chain $X(t)$ is called null-ergodic, if all
$p_i(t) \to 0$ $t \to \infty$ for any initial condition, and it is
called  weakly ergodic, if $\|{\bf p}^*(t)-{\bf p}^{**}(t)\| \to 0$
as $t \to \infty$ for any initial condition ${\bf p}^*(0), {\bf
p}^{**}(0)$.

\section{Bounds on the rate of convergence for a differential equation}

Consider  a general (linear or non-linear) differential equation
\begin{equation} \label{ur0100}
\frac{d\vy}{dt}=H\vy(t),
\end{equation}
in the space of sequences $l_1$ under the assumption of existence
and uniqueness of solution for any initial condition $\vy(0)$.

\smallskip

Let $H=(h_{ij})$, where all $h_{ij}$ depend on some parameters
 (for instance, on $y$, $t,\ldots$).

We have
$$\frac{d y_i}{dt}=h_{ii}y_i+\sum_{j\neq i} h_{ij}y_j.$$
Now, if $y_i > 0$, then
$$\frac{d |y_i|}{dt}= \frac{d y_i}{dt} = h_{ii}|y_i|+\sum_{j\neq i} h_{ij}y_j \le h_{ii}|y_i|+\sum_{j\neq i} |h_{ij}||y_j|,$$
and if $y_i < 0$, then we also have
$$\frac{d |y_i|}{dt}= - \frac{d y_i}{dt} = - h_{ii}y_i -\sum_{j\neq i} h_{ij}y_j \le h_{ii}|y_i|+\sum_{j\neq i} |h_{ij}||y_j|.$$

\smallskip

Finally, using the continuity of all coordinates of the solution and
the absolute convergence of all series, we obtain the estimate
\begin{eqnarray}
\frac{d \|y\|}{dt} = \sum_i \frac{d |y_i|}{dt} \le \sum_i\left(
h_{ii}|y_i|+\sum_{j\neq i} |h_{ij}||y_j| \right) \le \beta^* \|y\|,
\label{lognorm1}
\end{eqnarray}
\noindent where
\begin{equation}
\beta^*= \sup_i\left(h_{ii}+\sum_{j\neq i} |h_{ji}|\right).
\label{lognorm2}
\end{equation}

\smallskip

\begin{remark}\hspace{-0.2cm}{\bf .} One can see that inequality (\ref{lognorm1}) implies the
bound
\begin{equation}
\|y(t)\| \le e^{\int_0^t \beta^* \, du} \|y(0)\|. \label{lognorm3}
\end{equation}
Moreover, if $H$ is bounded for any $t$ linear operator function
from $l_1$ to itself, then $\beta^*(t)=\gamma(H(t))$ is the
corresponding logarithmic norm of $H(t)$, see
\cite{gz04,z95spa,z95jap,z06}.

 On the other hand, in a non-linear situation,
$\beta^*(t)$ yields a generalization of this notion.
\end{remark}

\section{Bounds on the rate of convergence for a projection of multidimensional BDP}

Again consider the forward Kolmogorov system (\ref{ur01}). Then we
have
\begin{eqnarray} \label{ur01*j}
\frac{dp_{{\bf m}}}{dt}= \sum_{l}\lambda _{l,{\bf
m-{\bf e_l}}}(t)p_{{\bf m}-{\bf e_l}} +\\
 \sum_{l}\mu_{l,{\bf m + {\bf e_l}}}(t)p_{{\bf m} + {\bf e_l}} -
\sum_{l}\left(\lambda _{l,{\bf m}} +\mu_{l,{\bf m}}
\right)(t)p_{{\bf m}}, \nonumber
\end{eqnarray}
\noindent for any ${\bf m}$.

\smallskip

In this section we consider the one-dimensional process $X_j(t)$ for
a fixed $j$. Denote $x_k(t)=\Pr\left(X_j(t)=k\right)$. Then
$x_k(t)=\sum_{{\bf m}, m_j=k}p_{{\bf m}}(t)$. The process $X_j(t)$
has nonzero jump rates only for unit jumps ($\pm 1$), namely, if
$X_j(t)=k$, then for small positive $h$ only the jumps
$X_j(t+h)=k\pm 1$ are possible with positive intensities, say
$\tilde{\lambda}_k$ and $\tilde{\mu}_k$ respectively. Moreover,
(\ref{ur01*j}) implies the equalities
\begin{equation} \label{ur01*j1}
\tilde{\lambda}_k x_k(t)= \sum_{{\bf m}, m_j=k}\lambda _{j,{\bf
m}}(t)p_{{\bf m}}(t),
\end{equation}
\begin{equation} \label{ur01*j2}
\tilde{\mu}_k x_k(t)= \sum_{{\bf m}, m_j=k}\mu_{j,{\bf
m}}(t)p_{{\bf m}}(t),
\end{equation}
\noindent and hence
\begin{equation}
\tilde{\lambda}_k=\frac{\sum_{{\bf m}, m_j=k}\lambda _{j,{\bf
m}}(t)p_{{\bf m}}(t)}{\sum_{{\bf m}, m_j=k}p_{{\bf m}}(t)},
\label{m01}
\end{equation}
\noindent  and
\begin{equation}
\tilde{\mu}_k=\frac{\sum_{{\bf m}, m_j=k}\mu _{j,{\bf m}}(t)p_{{\bf
m}}(t)}{\sum_{{\bf m}, m_j=k}p_{{\bf m}(t)}}. \label{m02}
\end{equation}

\smallskip

Then $X_j(t)$ is  an (in general, non-Markovian) birth and death
process with birth and death intensities $\tilde{\lambda}_k$ and
$\tilde{\mu}_k$ respectively (which depend on $t$ and the initial
condition of the original multidimensional process $X(t)$.)

\smallskip

For any fixed initial distribution ${\bf p}(0)$ and any $t>0$ the probability distribution ${\bf p}%
(t)$ is unique. Hence, $\tilde{\lambda}_k=\lambda _k\left( {\bf
p}(0),t\right) $ and $\tilde{\mu}_k=\mu _k\left( {\bf
p}(0),t\right)$ uniquely define the system
\begin{equation} \label{ur0101}
\frac{d\vx}{dt}=\tilde{A}\vx(t),
\end{equation}
for the vector $\vx(t)$ of state probabilities of the projection
$X_j(t)$ under the given initial condition. Here $\tilde{A}$ is the
corresponding three-diagonal ``birth-death'' transposed intensity
matrix with nonnegative for any $t$ and any initial condition ${\bf
p}(0)$ off-diagonal elements and zero column sums.

\smallskip

Let for all ${\bf m}$ and any $t \ge 0$
\begin{equation}
l_j\leq \lambda _{j,{\bf m}}(t)\leq L_j, \quad m_j\leq \mu _{j,{\bf
m}}(t)\leq M_j. \label{m03}
\end{equation}

Then from (\ref{m01}) and (\ref{m02}) we obtain the two-sided bounds
\begin{equation}
l_j\leq \tilde{\lambda}_k \leq L_j, \quad m_j\leq \tilde{\mu}_k \leq
M_j, \label{m03*}
\end{equation}
\noindent for any $k$, any $t$ and any initial conditions.

\smallskip

{\bf 1.}  Let
\begin{equation}
 M_j < l_j. \label{m08}
\end{equation}

Put $\sigma =\sqrt{M_j/l_j}<1$,  $\delta_n=\sigma ^n,\quad n\geq 0,$
 $\tilde{x}_n=\delta_n x_n$, and
 $\tilde{\vx}=\left(\tilde{x}_0,\tilde{x}_1,\dots\right)$. Let
 $\Lambda$ be a diagonal matrix,   $\Lambda = diag\left(\delta_0, \delta_1, \dots
 \right)$.

\smallskip

Then
\begin{equation} \label{ur0102}
\frac{d\tilde{\vx}}{dt}=\Lambda \tilde{A} \Lambda^{-1}
\tilde{\vx}(t).
\end{equation}
\noindent Hence
\begin{eqnarray}
\tilde{\lambda}_k+\tilde{\mu}_k-\frac{\delta_{k+1}}{\delta_k}\tilde{\lambda}
_k-\frac{\delta_{k-1}}{\delta_k}\tilde{\mu}_k\ge
\tilde{\lambda}_k\left( 1-\sigma \right) -\tilde{\mu}_k\left(
1/\sigma -1\right) \geq \\\nonumber l_j \left( 1-\sigma \right) -M_j
\left( 1/\sigma -1\right) = \left( \sqrt{l_j }-\sqrt{M_j }\right) ^2
= \alpha^{*}, \label{eq8}
\end{eqnarray}
\noindent we obtain the estimate
\begin{eqnarray}
\frac{d \|\tilde{\vx}\|}{dt} \le \sup_k \left(
\frac{\delta_{k+1}}{\delta_k}\tilde{\lambda}
_k+\frac{\delta_{k-1}}{\delta_k}\tilde{\mu}_k -
\tilde{\lambda}_k-\tilde{\mu}_k \right) = \nonumber \\ -\inf_k
\left(
\tilde{\lambda}_k+\tilde{\mu}_k-\frac{\delta_{k+1}}{\delta_k}\tilde{\lambda}
_k-\frac{\delta_{k-1}}{\delta_k}\tilde{\mu}_k \right) \le - \alpha^*
\|\tilde{\vx}\|, \label{lognorm1*}
\end{eqnarray}
\noindent and the following statement.

\smallskip

\begin{theorem}\hspace{-0.2cm}{\bf .} Let (\ref{m08}) hold for some $j$. Then
$X_j(t)$ is null-ergodic and the following bounds hold:
\begin{equation}
\|\tilde{\vx}(t)\| \le e^{ -\alpha^*  t} \|\tilde{\vx}(0)\|,
\label{lognorm3*}
\end{equation}
\noindent and
\begin{equation}
\Pr \left( X_j(t)\leq n/X_j(0)=k\right) \leq \sigma ^{k-n}\cdot
e^{-\alpha^*t}.  \label{m09}
\end{equation}
\end{theorem}

\smallskip
{\bf 2.}  Let
\begin{equation}
L_j < m_j, \quad \alpha_*= l_j+m_j -2\sqrt{L_jM_j} > 0. \label{m04}
\end{equation}

\smallskip

The property ${\bf x}(t) \in \Omega$ for any $t \ge 0$ allows to set
$x_0(t) = 1 - \sum_{i \ge 1} x_i(t)$. Then from (\ref{ur0101}) we
obtain the system
\begin{equation}\label{216}
\frac{d{\bf z}}{dt}=\tilde{B}{\bf z}+\tilde{{\bf f}},
\end{equation}
\noindent  where  ${\bf z} = (x_1,x_2,\dots)^\top$, $\tilde{{\bf f}}
= \big(\tilde{\lambda}_0,0,0,\dots\big)^\top$, and the corresponding
matrix $\tilde{B} = \big(\tilde{b}_{ij}\big)_{i,j=1}^{\infty}$, and
$\tilde{b}_{ij} =\tilde{a}_{ij}-\tilde{a}_{i0}$ for the
corresponding elements of the matrix $\tilde{A}$.

\smallskip

For the solutions of system (\ref{216}) the rate of convergence is
determined by the system
\begin{equation}\label{216o}
\frac{d{\bf w}}{dt}=\tilde{B}{\bf w},
\end{equation}
\noindent where all elements of $\tilde{B}$ depend on $t$ and
initial condition of the original process.
\smallskip

Now let $\beta= \sqrt{\frac{M_j}{L_j}}>1$ in accordance with
(\ref{m04}). Let $d_{k+1}=\beta^k$, $k \ge 0$. Denote  by $D$ the
upper triangular matrix
\begin{equation}
D=\left(
\begin{array}{ccccccc}
d_1   & d_1 & d_1 & \cdots  \\
0   & d_2  & d_2  &   \cdots  \\
0   & 0  & d_3  &   \cdots  \\
& \ddots & \ddots & \ddots \\
\end{array}
\right). \label{204}
\end{equation}
\noindent Let  $\tilde{\bf w}= D{\bf w}$. Then the following bound
holds:
\begin{eqnarray}
\frac{d \|\tilde{\bf w}\|}{dt} \le \sup\limits_{i \ge
0}\left(\frac{d_{i+1}}{d_i}
 \tilde{\lambda}_{i+1} + \frac{d_{i-1}}{d_i}  \tilde{\mu}_i - \left( \tilde{\lambda}_i + \tilde{\mu}_{i+1}
\right)   \right) \nonumber  =
\\  - \inf \limits_{i \ge 0} \left(\left( \tilde{\lambda}_i + \tilde{\mu}_{i+1} - \beta \tilde{\lambda}_{i+1} -\tilde{\mu}_{i}/\beta \right) \right)  \le -
\alpha_*\|\tilde{\bf w}\|,  \label{219}
\end{eqnarray}
\noindent and we obtain the following statement.

\smallskip

\begin{theorem}\hspace{-0.2cm}{\bf .} Let (\ref{m04}) hold for some $j$. Then
$X_j(t)$ is weakly ergodic and the following bound holds:
\begin{equation}
\|D{\bf w}(t)\| \le e^{- \alpha_* t}\|D{\bf w}(0)\|, \label{214'}
\end{equation}
\noindent for any $ t  \ge 0$ and any corresponding initial
conditions.
\end{theorem}

\smallskip

\begin{remark}\hspace{-0.2cm}{\bf .} Instead of $X_j(t)$ we can obtain the same results
for the one-dimensional process $Z(t) = |X(t)|$, that is, the number
of all  particles at the moment $t$.

\end{remark}

\section{ACKNOWLEDGMENTS}

The work of Zeifman and Korolev is supported by the Russian Science Foundation under grant 18-11-00155,
they obtained bounds on the rate of convergence in weakly ergodic situation.

\end{document}